\documentclass[11pt]{article}
\usepackage{amsmath, amsthm, amssymb}    
\usepackage{bridges}                     
\usepackage{graphicx}                    
\usepackage{hyperref}                    
\usepackage{enumitem}

\usepackage{comment} 

\usepackage{xcolor}

\title{A Gallery of Gaussian Periods}

\author{Ellen Eischen\textsuperscript{1} and Stephan Ramon Garcia\textsuperscript{2}
\vspace{10pt}\\
\textsuperscript{1}Department of Mathematics, University of Oregon,
{\tt eeischen@uoregon.edu}\\
\textsuperscript{2}Department of Mathematics, Pomona College,
{\tt stephan.garcia@pomona.edu}
}

\date{}

\usepackage{geometry}
\usepackage{graphicx}
\usepackage{float}
\usepackage{hyperref}
\usepackage{ellipsis}
\usepackage{amsmath,amsthm}
\usepackage{enumerate}
\usepackage{subcaption}
\usepackage{url}
\usepackage{cite}

\allowdisplaybreaks

\graphicspath{{./pics/}}

\newcommand{\G}{\mathsf{G}}
\newcommand{\Q}{\mathbb{Q}}

\newcommand{\C}{\mathbb{C}}
\newcommand{\T}{\mathbb{T}}
\newcommand{\Z}{\mathbb{Z}}

\renewcommand{\pmod}[1]{\,(\operatorname{mod}#1)}

\renewcommand{\phi}{\varphi}

\newcommand{\Aut}{\operatorname{Aut}}

\newcommand{\IQ}{\Q}
\newcommand{\Gal}{\mathrm{Gal}}
\newcommand{\ZZ}{\Z}
\newcommand{\isomto}{\overset{\sim}{\rightarrow}}

\newcommand{\GPset}{\G(n, \omega)}
\newcommand{\CO}{\mathcal{O}}
\newcommand{\ord}{\mathrm{ord}}

\theoremstyle{definition}

\newtheorem{Example}{Example} 

\allowdisplaybreaks

\begin{document}

\maketitle
\thispagestyle{empty}

\begin{abstract}
Gaussian periods are certain sums of roots of unity whose study dates back to Gauss's seminal work in algebra and number theory.  Recently, large scale plots of Gaussian periods have been revealed to exhibit striking visual patterns, some of which have been explored in the second named author's prior work.  In 2020, the first named author produced a new app, \texttt{Gaussian periods}, which allows anyone to create these plots much more efficiently and at a larger scale than before.  In this paper, we introduce Gaussian periods, present illustrations created with the new app, and summarize how mathematics controls some visual features, including colorings left unexplained in earlier work.
\end{abstract}
\bigskip

\section*{Introduction}

Gaussian periods, certain sums of roots of unity introduced by Gauss, have played a key role in several mathematical developments.  For example, Gauss employed them in his work on constructibility of regular polygons by unmarked straightedge and compass, as well as in number theory.  In the past few years, realizations of their role in the supercharacter theory of Diaconis and Isaacs have led in new directions \cite{SESUP}.  

\begin{figure}[h!tbp]
\centering
\begin{minipage}[b]{0.475\textwidth} 
	\includegraphics[width=\textwidth]{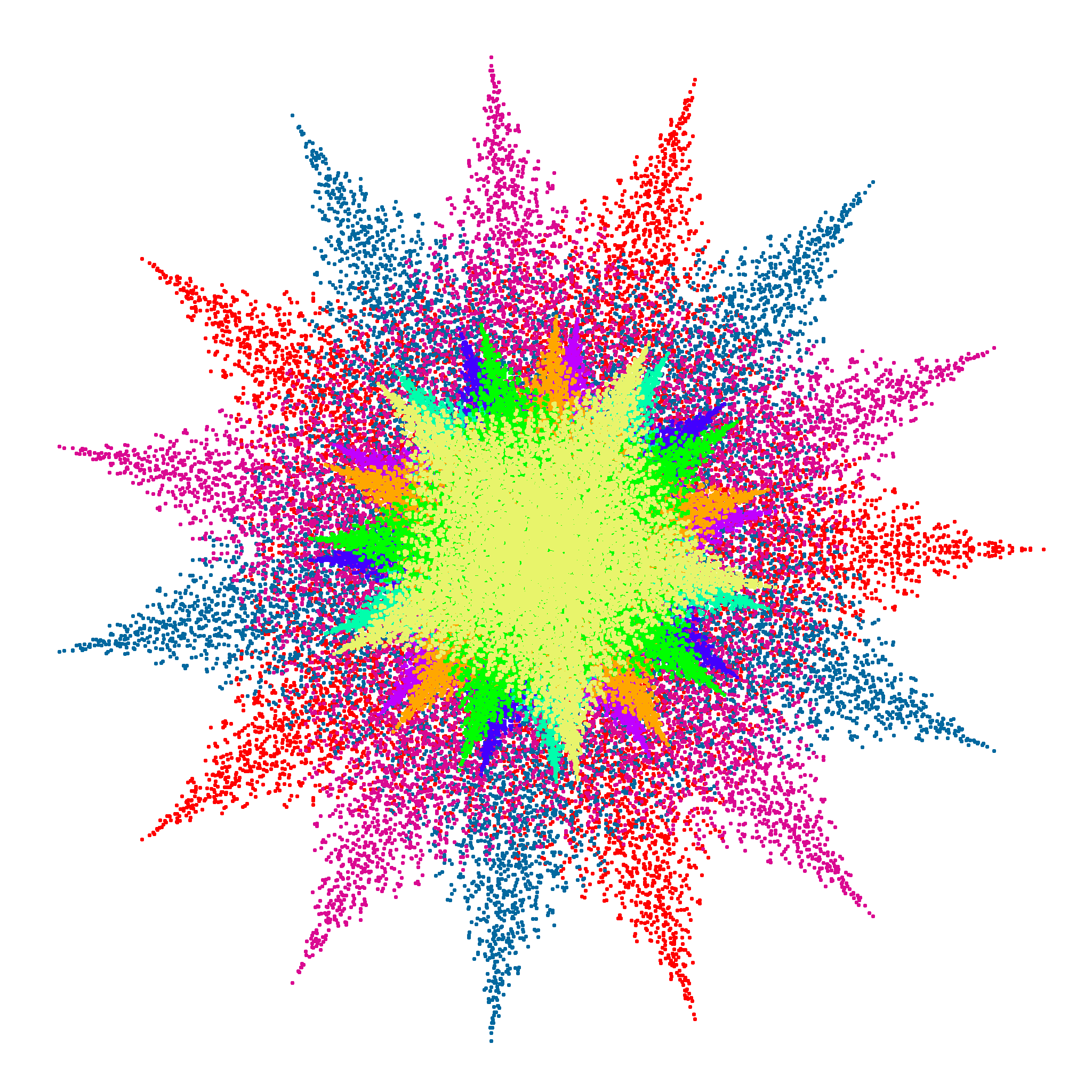}
        	\subcaption{$n=1481151 $, $\omega = 54184$, $c=21$} 
\end{minipage}
\hfill
\begin{minipage}[b]{0.475\textwidth} 
		\includegraphics[width=\textwidth]{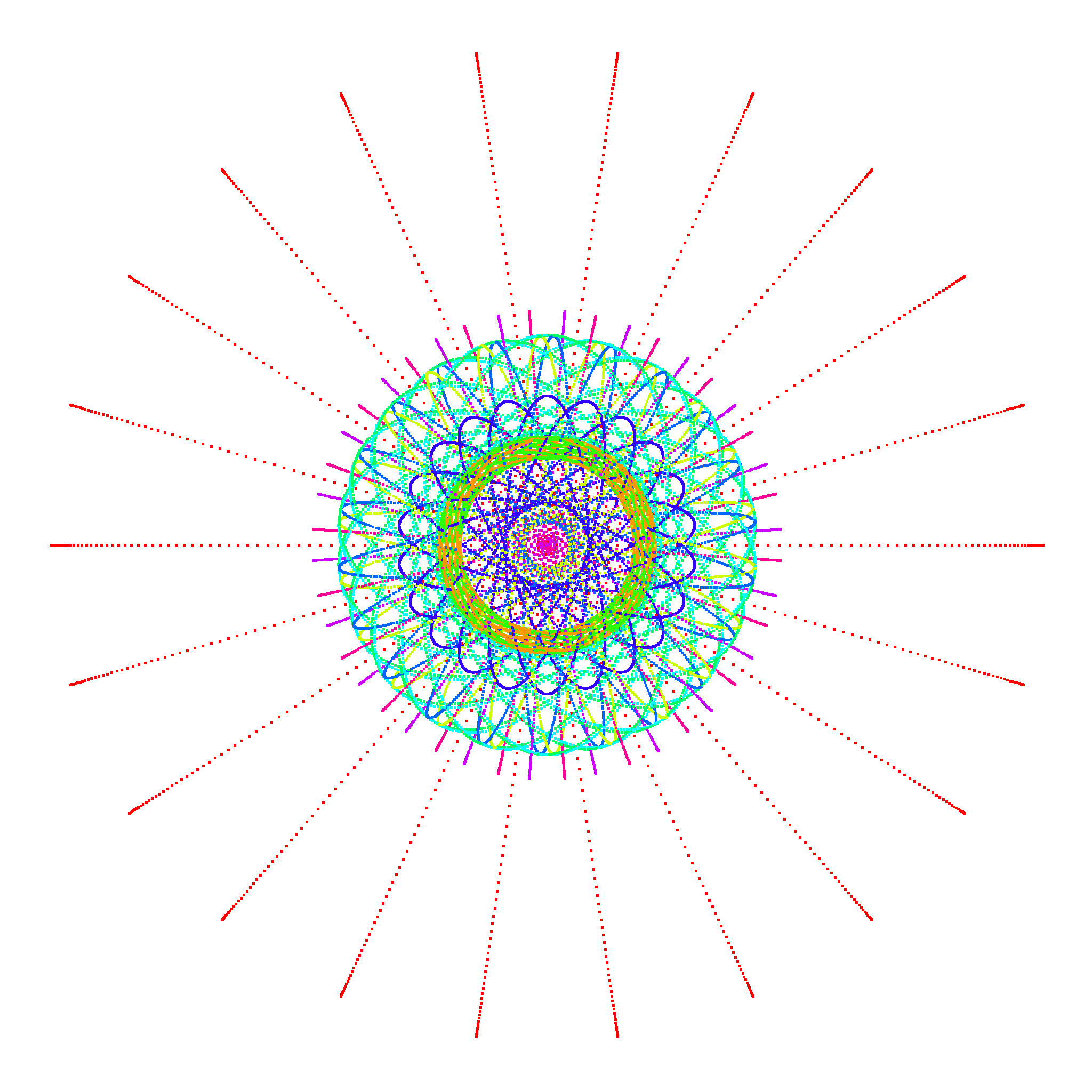}
        	\subcaption{$n=255255$, $\omega = 254$,  $c=7$} 
        	\label{Figure:Rotate:b}
\end{minipage}
\caption{Examples of $\GPset$ for different choices of input.}
\label{Figure:BeautifulGallery}
\end{figure}

The fast computations afforded by modern technology provide new insights into Gaussian periods by enabling us to study them at a scale that was until recently unfathomable.  In particular, large-scale plots of Gaussian periods display striking visual properties that were previously unknown (see, e.g., Figure \ref{Figure:BeautifulGallery}).
The goals of this paper lie in the intersection of the mathematical and artistic aspects of Gaussian periods:
\begin{itemize}[leftmargin=*]
\item{Introduce Gaussian periods, together with new illustrations representing some of their visual features.}
\item{Introduce a new app, created in 2020, that quickly produces illustrations of Gaussian periods and is appropriate as a tool for art, exploratory mathematical research, illustration at scale, and pedagogy.}
\item{Summarize how mathematics controls some features, including colorings unexplained in earlier work.}\label{item:coloringsadhoc}
\end{itemize}

Before proceeding, we clarify what we mean by the term {\em Gaussian period}.  Given a positive integer $n$, an integer $\omega$ coprime to $n$, and an integer $k$, we set
\begin{align*}
\eta_{n, \omega, k}:=\sum_{j=0}^{d-1}e^{\frac{2\pi i\omega^j k}{n}},
\end{align*}
where $d$ is the multiplicative order of $\omega$ mod $n$; that is, $d$ is the smallest positive integer such that $\omega^d\equiv 1\pmod{n}$.  When $k$ is relatively prime to $n$, we say that $\eta_{n, \omega, k}$ is a {\em Gaussian period of modulus $n$ and generator $\omega$}.
In this paper, by a slight abuse of notation that coincides with the conventions employed in \cite{GHM,GNGP, Lutz}, we still call $\eta_{n, \omega, k}$ a {\em Gaussian period of modulus $n$ and generator $\omega$} even when $k$ is not coprime to $n$.  As we explain later in this paper, $\eta_{n, \omega, k}$ is a positive integer multiple of the Gaussian period $\eta_{\frac{n}{\gcd(n, k)}, \omega, k}$, which allows us to preserve certain key structures of interest.

\section*{Creating and viewing plots of Gaussian periods}
Each of the images in this paper is a plot in the complex plane of
\begin{align*}
\G(n, \omega) := \left\{\eta_{n, \omega, k} : k = 1, 2, \ldots, n\right\}\subset\mathbb{C}.
\end{align*}
Mathematical principles guarantee that $\G(n, \omega)$ exhibits certain basic symmetries.  For large values of $n$ (which we refer to as {\em large scale}), however, plots of $\G(n, \omega)$ also exhibit striking patterns and intricacies whose aesthetic properties can be appreciated even by those without mathematical training.  

That large-scale plots of Gaussian periods exhibit such variety and intricate patterns was a surprise.  This was discovered by B.\ Lutz,
an undergraduate, in the course of his senior thesis.  In his explorations, he also introduced a coloring scheme, which is employed in the plots in this paper and discussed in the next section.  

A new app, \texttt{Gaussian Periods} (written in Swift for Apple computers and freely available \cite{app}) was produced in 2020 by the first named author, with assistance from R.\ Lipshitz.  This app:
\begin{itemize}[leftmargin=*]
\item{plots Gaussian periods faster than previous code (including the aforementioned \texttt{Mathematica} code), taking seconds to produce images that used to take hours;}
\item{allows larger scale plots than previously possible (e.g., Figure \ref{Figure:Fills:c}, which contains over 9 million points), which can be useful for exploring or illustrating asymptotic behavior;}
\item{does not require programming experience or mathematical expertise;}
\item{allows one to quickly modify values of $n$ and $\omega$, as well as a coloring parameter $c$; and}
\item{includes an option to save layers suitable for further steps, e.g., manipulation in \texttt{Adobe Photoshop} to customize color choices (as was done to produce the images in this paper).}
\end{itemize}
As a result, the app is suitable for projects ranging from art to exploratory mathematical investigations to illustration.  We also used it to produce all the images in this paper.  To improve the image quality in this paper, we layered the different color components (plotting all the points corresponding to a color at the same time), which is assisted by the layers option in the app.

  While notions of ``beauty'' and ``aesthetic appeal'' are subjective, symmetry has often appeared in discussions of beauty in both art and math, dating back to the ancient Babylonians.   In a historical context, given the algebraic origins of Gaussian periods, it is fitting to note that symmetry originally became a prominent concept in mathematics not through geometry, but through algebra \cite[Preface]{stewart}.  Gaussian periods, together with the various symmetries they exhibit, can also be considered of independent artistic merit and be appreciated in their own right (even by those with no mathematical training).  More broadly, patterns in plots of certain other families of algebraic numbers have also been recognized for their beauty \cite{baez}.  
  For those who wish to delve further into philosophical considerations of the artistic merits of Gaussian periods, a brief survey of visual aesthetics in similar mathematical contexts can be found in \cite[pp.~121-4]{Brown}. 
  
  \section*{Using color to reveal structures}
Whether focusing on art or math, a coloring scheme can be used to highlight some structures in $\G(n, \omega)$ for a given $n$ and $\omega$.  Following \cite[\S 3]{GHM}, 
fix a positive integer $c \mid n$ and, for $j=1,2, \ldots, c$, assign the same color to all points in the set
$\{\eta_{n, \omega, k} : k\equiv j\pmod{c}\}$.
Two points $\eta_{n, \omega, k}$ and $\eta_{n, \omega, \ell}$ might have the same color even if $k\nequiv \ell\pmod{c}$, since $\eta_{n, \omega, k} = \eta_{n, \omega, k\omega^j}$ for all integers $j$. 

Since this is the coloring scheme employed in each of the images here (and is implemented in \texttt{Gaussian Periods}), we begin with small-scale examples to help readers grasp the subtleties of this coloring scheme. 

\begin{figure}[h!tbp]
\centering
\begin{minipage}[b]{0.3\textwidth} 
	\boxed{\includegraphics[width=\textwidth]{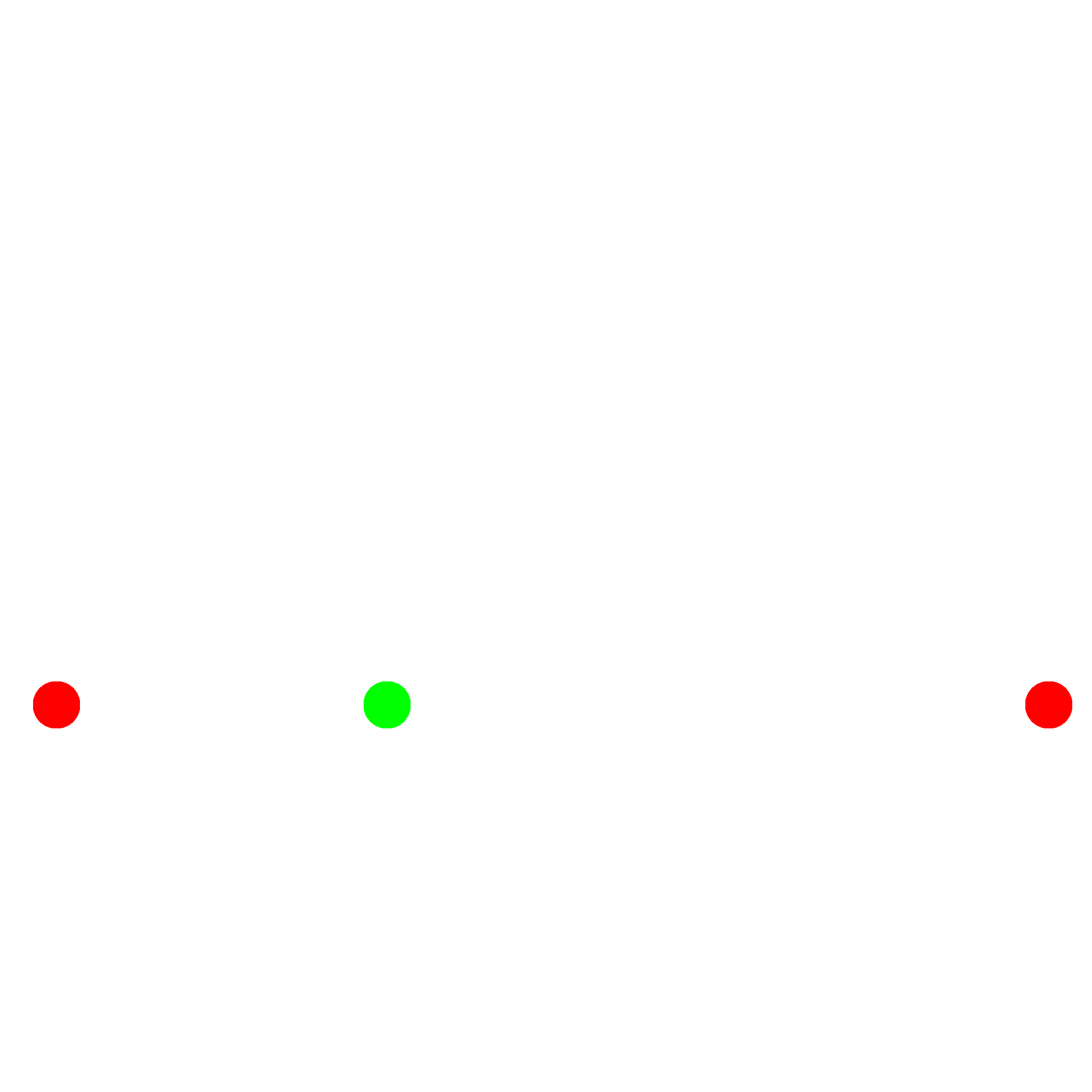}}
        	\subcaption{$n=27$, $\omega = 2$, $c=9$} 
        	\label{Figure:Small:a}
\end{minipage}
\hfill
\begin{minipage}[b]{0.3\textwidth} 
	\boxed{\includegraphics[width=\textwidth]{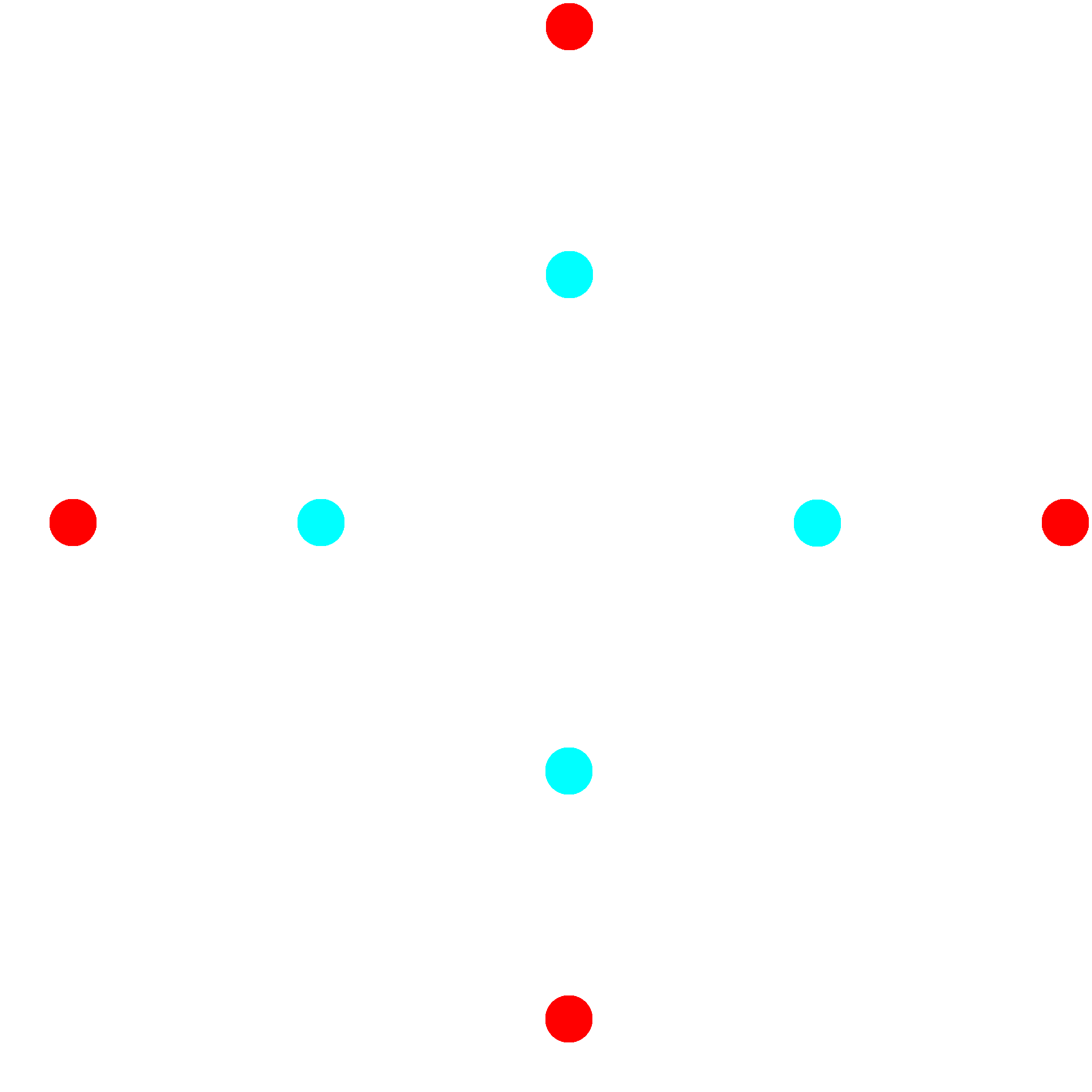}}
        	\subcaption{$n=12$, $\omega = 5$, $c=3$} 
        	\label{Figure:Small:b1}
\end{minipage}
\hfill
\begin{minipage}[b]{0.3\textwidth} 
	\boxed{\includegraphics[width=\textwidth]{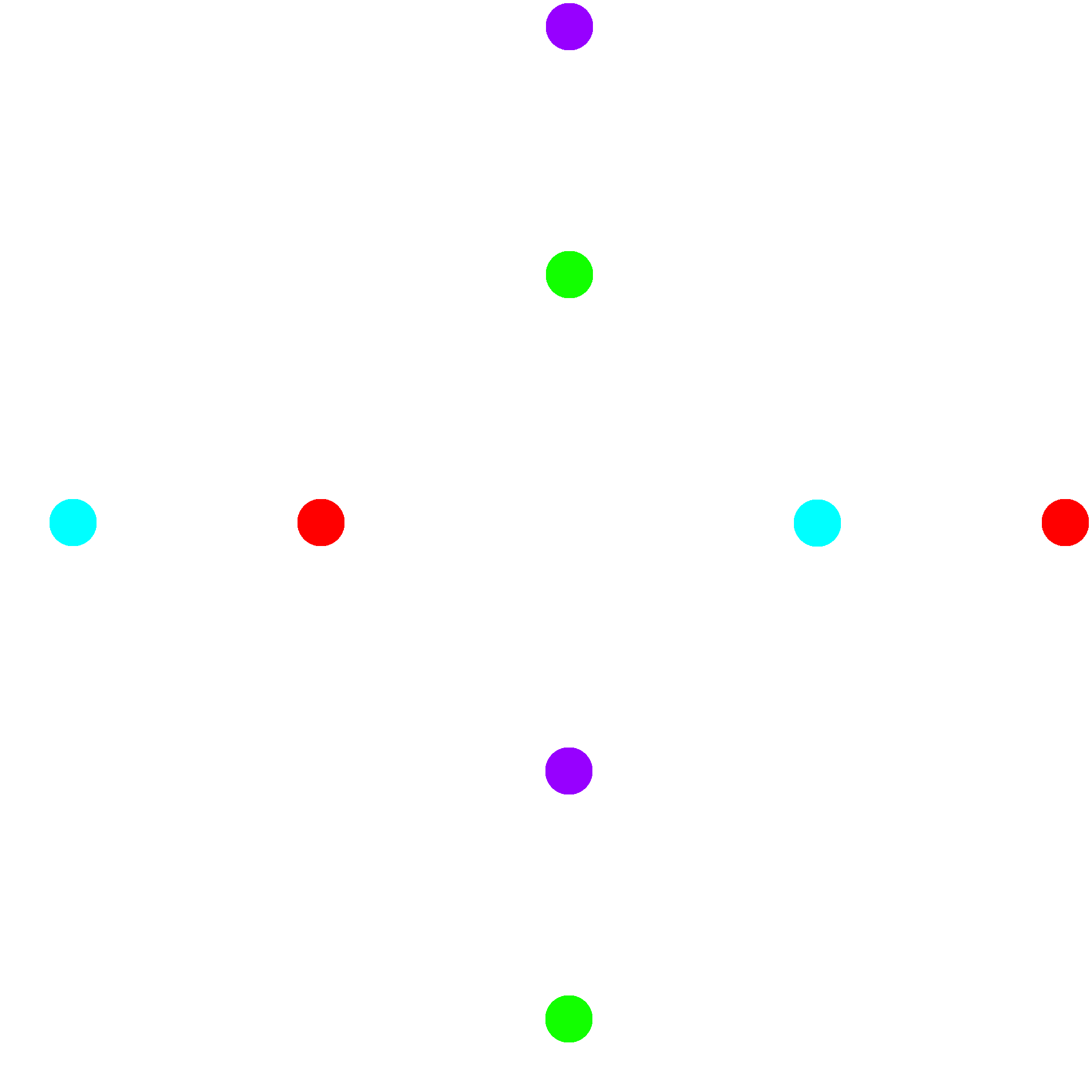}}
        	\subcaption{$n=12$, $\omega = 5$, $c=4$} 
        	\label{Figure:Small:b2}
\end{minipage}
\caption{Small scale illustrations for Examples \ref{example:small:a} and \ref{example:small:b}.}
\label{Figure:Small}
\end{figure}

\begin{Example}\label{example:small:a}
Suppose $n = 27$, $\omega = 2$, and $c = 9$.  There are three orbits of $\langle\omega\rangle$ acting on $\ZZ/27\ZZ$: $\CO_1:=\ZZ/27\ZZ^\times$, $\CO_0:=\left\{0\right\}$, and $\CO_{18}:=\{[a]\in\ZZ/27\ZZ :  27>\gcd(a, 27) >1\}$, and $G(27, 2)$ consists of the three points $\eta_{27, 2, 1},$ $\eta_{27, 2, 0}$, and $\eta_{27, 2, 18}$.  So by our rules for coloring, since $c=9$ and $18\equiv 0\pmod{9}$, we must assign the same colors to $\eta_{27, 2, 0}$ and $\eta_{27, 2, 9}$.  Since elements in $(\ZZ/27\ZZ)^\times$ are not congruent modulo $9$ to elements divisible by $3$, there are no restrictions on the color of $\eta_{27, 2, 1}$.  As seen in Figure \ref{Figure:Small:a}, for this input, our recipe produces a plot with $3$ points and just $2$ colors, even though $c= 9$. 
\end{Example}

This approach to coloring was chosen not as a way to illustrate a particular mathematical principle but rather as a recipe Lutz discovered that produced visually appealing pictures.  In fact, in the final paragraph of \cite{GHM}, he and his coauthors note that this approach to coloring is \emph{ad hoc} and that ``a general theory is necessary to formalize our intuition.''

\section*{Algebraic structures}

Through the lens of Galois theory  (which was never considered in prior work on illustrating $\G(n,\omega)$), the above coloring scheme can be linked to mathematical structure.  The main theorem of Galois theory gives a bijection between the subfields of $\IQ\left(\zeta_n\right)$, with $\zeta_n$ a primitive $n$th of unity, and the subgroups of the Galois group $G:=\Gal\left(\IQ\left(\zeta_n\right)/\IQ\right)=\Aut\left(\IQ\left(\zeta_n\right)\right)$, which is identified with $\left(\ZZ/n\ZZ\right)^\times$ via
\begin{align*}
\left(\ZZ/n\ZZ\right)^\times&\isomto G\\
[m]&\mapsto \psi_m\\
\psi_m\Big(\sum_j a_j \zeta_n^j\Big) &:= \sum_j a_j\left(\zeta_n^m\right)^j, \hspace{0.5in} a_j\in\IQ.
\end{align*}
Each subgroup $H\subseteq G$ corresponds to the fixed field $\IQ\left(\zeta_n\right)^H:=\left\{z\in  \IQ\left(\zeta_n\right) : \psi(z) = z \mbox{ for all } \psi\in H\right\}$.  Then $H$ is the Galois group of $\IQ\left(\zeta_n\right)$ over $\Q\left(\zeta_n\right)^H$, i.e., the group of field automorphisms of $\IQ\left(\zeta_n\right)$ fixing $\IQ\left(\zeta_n\right)^H$, and $G/H$ is the Galois group of $\IQ\left(\zeta_n\right)^H$ over $\IQ$, i.e., the group of field automorphisms of $\IQ\left(\zeta_n\right)^H$. 

To help identify elements of the subfields of $\IQ\left(\zeta_n\right)$, students in Galois theory courses are sometimes assigned to determine Gaussian periods (e.g., see \cite[\S 14.5]{DF}).  Indeed,  if $H=\langle\omega\rangle$, then
\begin{equation*}
 \eta_{n, \omega, k}=\sum_{\psi\in H}{\left(\psi\left(\zeta_n\right)^k\right) }= \sum_{\psi\in H}\psi\big(\zeta_n^k\big)\in \IQ\left(\zeta_n\right)^H,
\end{equation*}
So 
the plot of $\G(n, \omega)$ contains a (rescaled) plot of $\G(\frac{n}{\gcd(n, k)}, \omega)$.  
If $c$ is the coloring number, then the plot of $\G(n, \omega)$ contains a single-colored plot of $\G(\frac{n}{c}, \omega)$, rescaled by a factor of $\frac{\ord_n(\omega)}{\ord_{n/c}(\omega)}$, with $\ord_m(\omega)$ denoting the multiplicative order of $\omega \pmod{m}$.   For clarification, we illustrate this in a simple example.

\begin{Example}\label{example:small:b}
Let $n= 12$ and $\omega = 5$.  Note that $5$ has multiplicative order $2$ mod $12$ and $1$ mod $4$.  So $\eta_{12, 5, 3k} = 2\eta_{4, 5, k} = 2 e^{\frac{2\pi ik}{4}}$ for all $k$.  So $G(4, 5)=\{1, -1, i, -i\}$, and $G(12, 5)\supseteq 2\cdot G(4, 5) = \{2, -2, 2i, -2i\}$.  If we choose $c = 3$, then the points in $2\cdot G(4, 5)$ all must be the same color as each other (the red, outer diamond in Figure \ref{Figure:Small:b1}), while the other points need not be colored that color (the blue, inner diamond).  On the other hand, selecting $c= 4$ forces the pair of points in $G(3, 5)\subset G(12, 5)$ to be the same color (shown in Figure \ref{Figure:Small:b2} in red).  As illustrated in Figure \ref{Figure:Small:b2}, the pair gets rotated by $2\pi/4 = 2\pi\omega/4$, with a new color allowed at each rotation.
\end{Example}

Regarding the current coloring scheme, if two elements of $\ZZ/n\ZZ^\times\cong\Gal(\IQ(\zeta_n)/\IQ)$ are congruent $\mod c$, then the corresponding elements of the Galois group restrict to the same element of $\Gal(\IQ(\zeta_c)/\IQ)$.  So, for example, given an element $\psi\in \Gal(\IQ(\zeta_c)/\IQ)$ and an integer $a$, all points $\Psi(\eta_{n, \omega, a})$, as $\Psi$ ranges over the extensions of $\psi$ to $\IQ(\zeta_n),$ are colored the same.  There are also other natural coloring schemes.  For example, an option (called ``period squared'') in \texttt{Gaussian Periods} is to color $\eta_{n, \omega, k}$ and $\eta_{n, \omega, -k}$ the same, thus creating symmetry across the real axis.  This corresponds to coloring the points in the Galois orbit of the complex conjugation automorphism the same.  More generally, one might color all points in some other given Galois orbits the same.  Furthermore, our reformulation in terms of elements of fixed fields should naturally generalize beyond Gaussian periods to the illustration of symmetries in other settings beyond the scope of this short paper (such as finite nonabelian extensions of $\IQ$).

\section*{Asymptotic behavior}
\begin{figure}[h!tbp] 
\centering
\begin{minipage}[b]{0.32\textwidth} 
	\includegraphics[width=\textwidth]{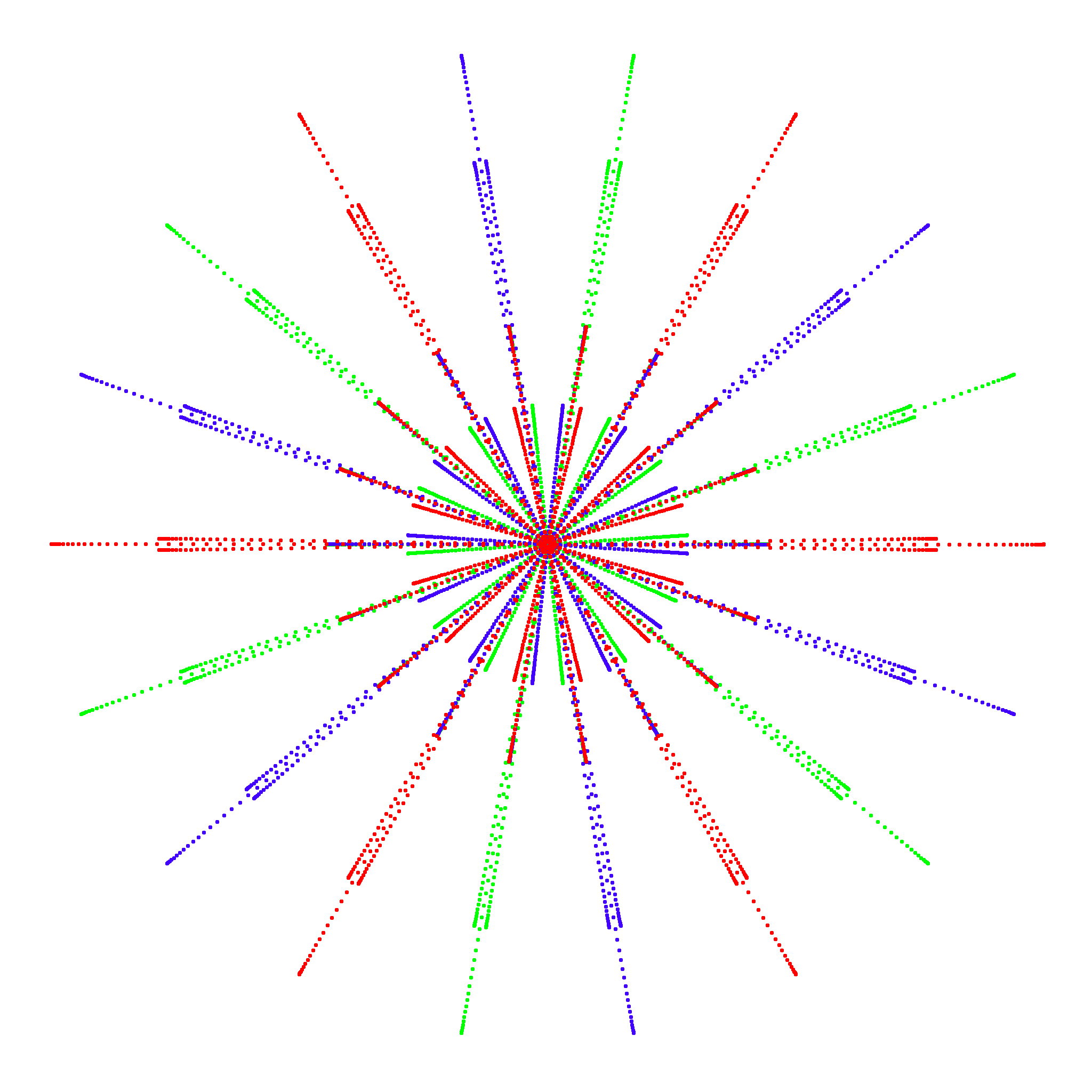}
        	\subcaption{$n=29070$, $\omega = 1189$, $\gcd(1188, 29070)=18$, $c=3$} 
        	\label{Figure:Rotate:a}
\end{minipage}
~ 
\begin{minipage}[b]{0.32\textwidth} 
	\includegraphics[width=\textwidth]{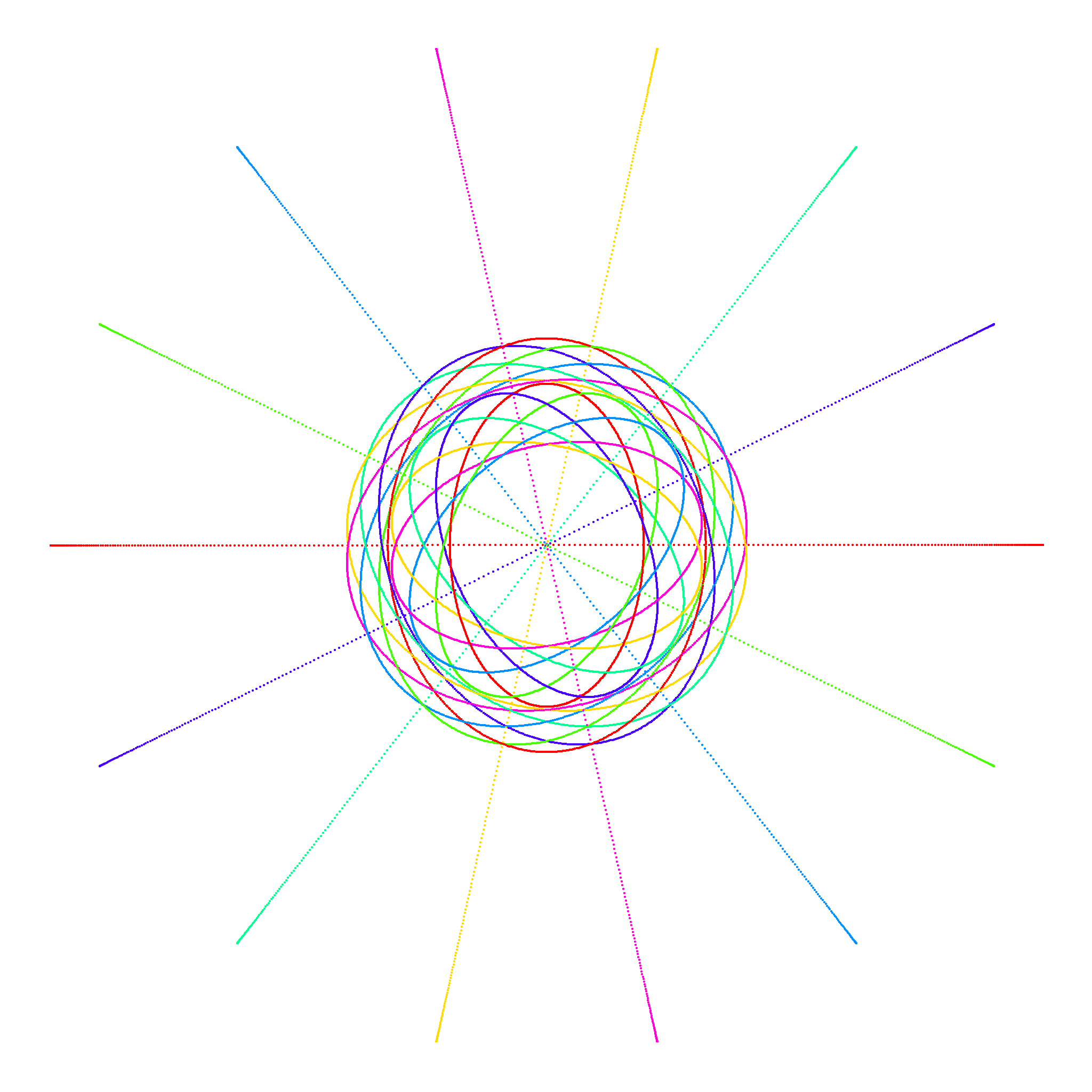}
        	\subcaption{$n=70091$, $\omega = 21792$, $\gcd(21791,70091)=7$, $c=7$} 
        	\label{Figure:Rotate:c}
\end{minipage}
~ 
\begin{minipage}[b]{0.32\textwidth} 
	\includegraphics[width=\textwidth]{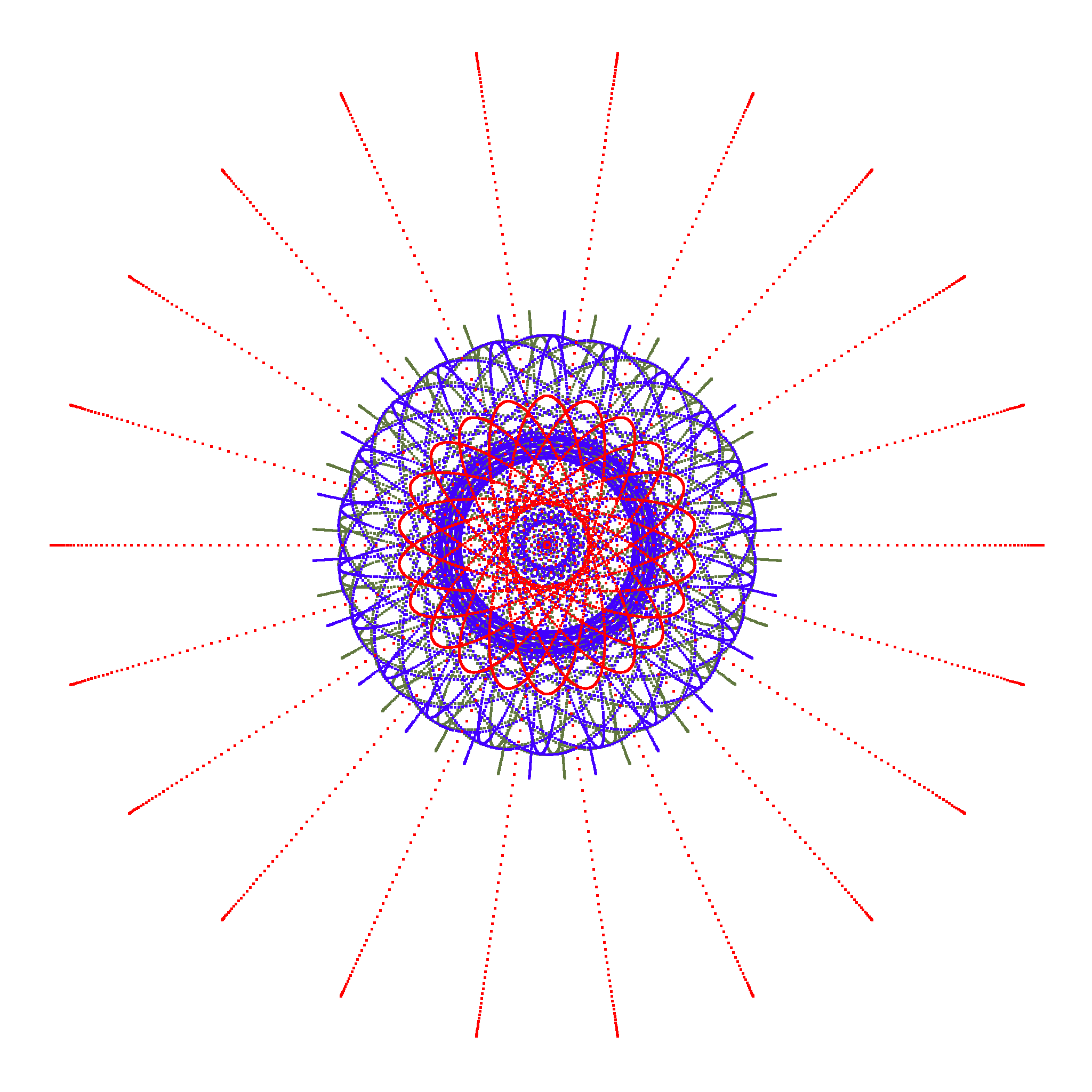}
        	\subcaption{$n=255255$, $\omega = 254$, $\gcd(255255,253)=11$, $c=7$} 
\end{minipage}
\caption{Dihedral symmetry of $\GPset$.}
\label{Figure:Rotate}
\end{figure}

Symmetry is the most obvious feature in typical Gaussian-period plots.
We say that $\GPset$ has \emph{$k$-fold dihedral symmetry} 
if $\GPset$ is invariant under the action of the dihedral group of order $2k$. 
That is, $\GPset$ is invariant under complex conjugation and rotation by $2\pi/k$ about the origin. 
It turns out that $\GPset$ has at least $\gcd(\omega-1,n)$-fold dihedral symmetry \cite[Prop.~3.1]{GNGP}.  
This symmetry refers to the uncolored graph; the colors highlight additional features beyond the initial symmetry.
This is illustrated in Figure \ref{Figure:Rotate}.

\begin{figure}[h!tbp]
\centering
\begin{minipage}[b]{0.32\textwidth} 
	\includegraphics[width=\textwidth]{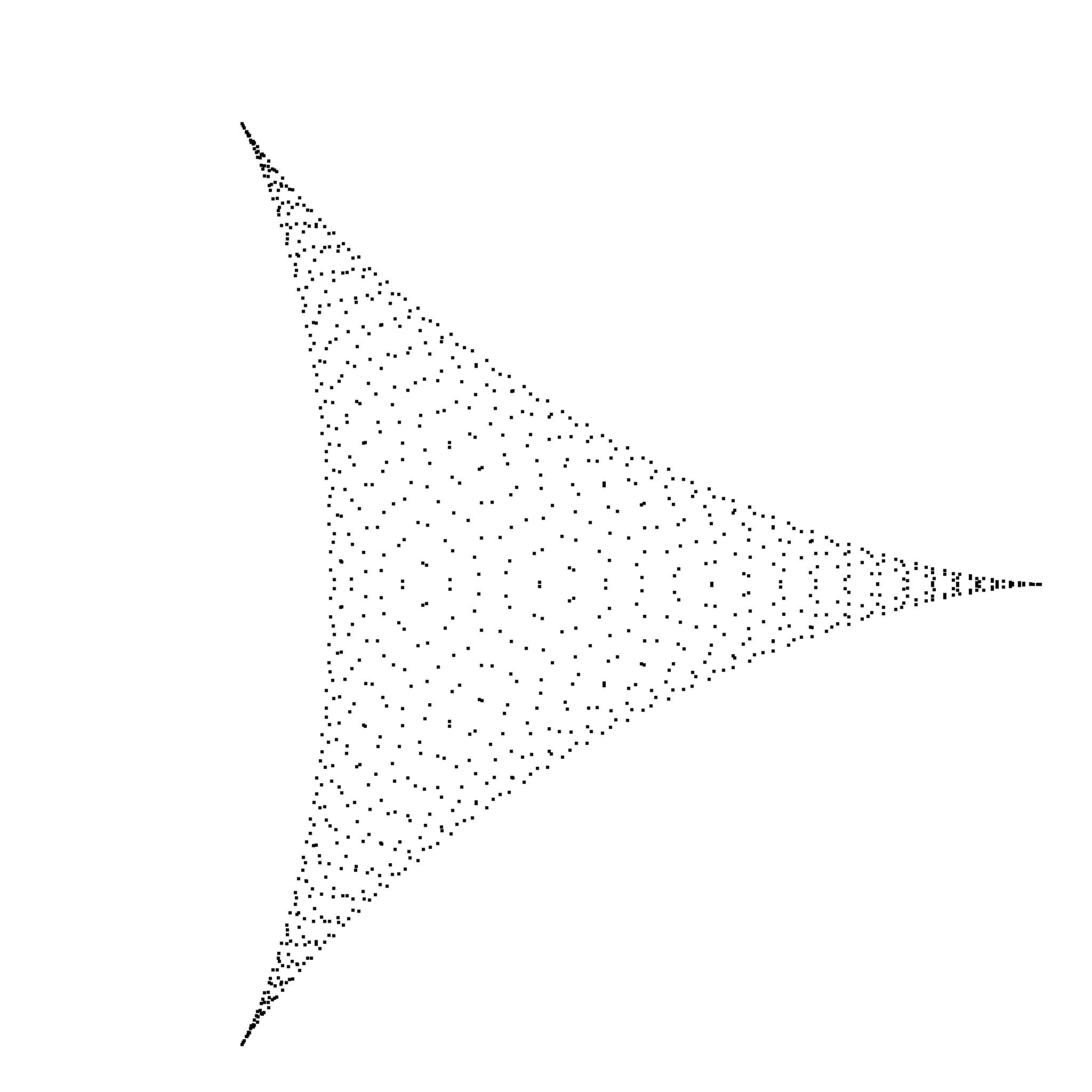}
        	\subcaption{$n=3019$, $\omega = 239$} 
\end{minipage}
~ 
\begin{minipage}[b]{0.32\textwidth} 
	\includegraphics[width=\textwidth]{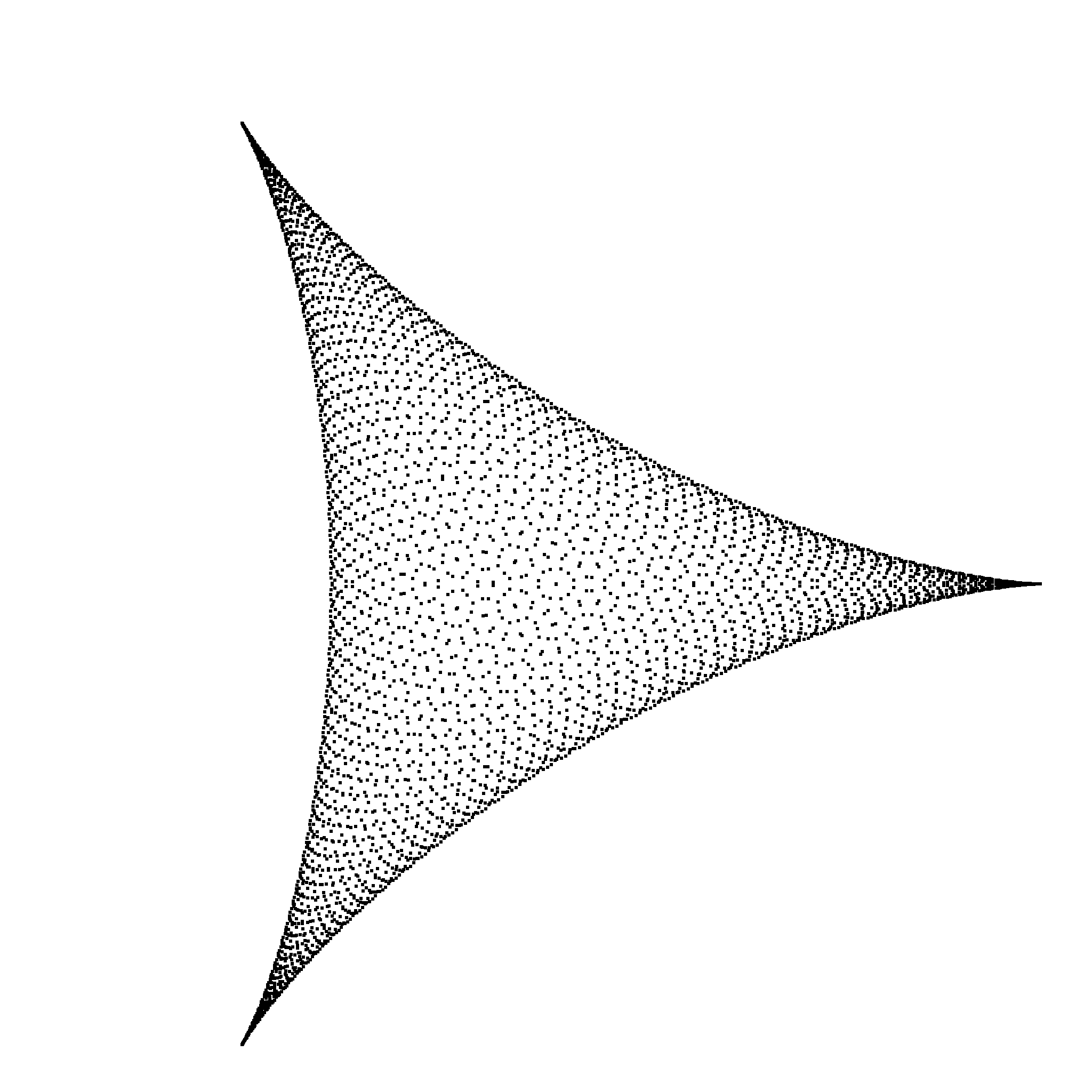}
        	\subcaption{$n=13063$, $\omega = 1347$} 
\end{minipage}
~ 
\begin{minipage}[b]{0.32\textwidth} 
	\includegraphics[width=\textwidth]{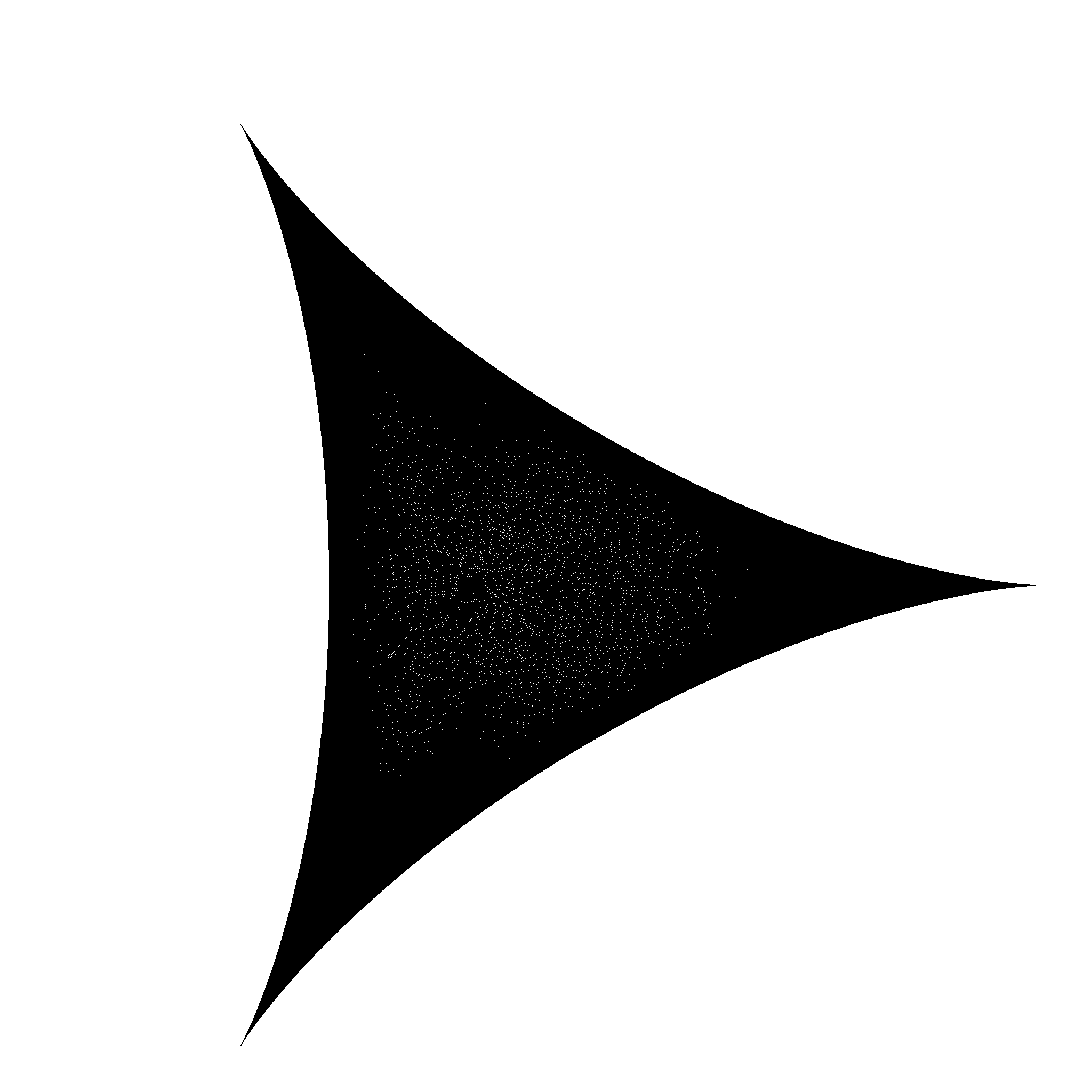}
        	\subcaption{$n=9114361$, $\omega = 3082638$}
	\label{Figure:Fills:c} 
\end{minipage}
\caption{Sometimes $\G(n,\omega)$ appears to ``fill out'' the image of a Laurent
polynomial $g:\T^{\phi(d)}\to\C$.  Here $g(z_1, z_2) = z_1 + z_1 + 1/(z_1 z_2)$.}
\label{Figure:Fills}
\end{figure}

Gaussian period plots often demonstrate great structural coherence if the parameters 
$n$ and $\omega$ vary in the appropriate manner.  The following ``filling out'' of various shapes
was discovered in \cite[Thm.~6.3]{GNGP}. See the thorough exposition in \cite[Thm.~1]{GHM}.
Let $q=p^a$ be a nonzero power of an odd prime and let $\omega = \omega(q)$ be such that
$d= |\langle \omega \rangle|$ divides $p-1$.  
Then $\G(q,\omega)$ is contained in the image of the Laurent polynomial function $g:\T^{\phi(d)}\to\C$ defined by
\begin{equation*}
	g(z_1,z_2,\ldots,z_{\phi(d)}) = \sum_{k=0}^{d-1} \prod_{j=0}^{\phi(d)-1} z_{j+1}^{b_{k,j}},
\end{equation*}
where the integers $b_{k,j}$ are determined by
	$t^k\equiv \sum_{j=0}^{\phi(d)-1} b_{k,j} t^j \pmod{\Phi_d(t)}.$
Here $\T$ denotes the unit circle in $\C$, $\phi$ is the Euler totient function, and $\Phi_d$ denotes the $d$th cyclotomic polynomial.
For a fixed $d$, as $q$ becomes large, $\G(q,\omega)$ ``fills out'' the image of $g$; see Figure \ref{Figure:Fills}.  

\begin{figure}[h!tbp]
\centering
\begin{minipage}[b]{0.32\textwidth} 
	\includegraphics[width=\textwidth]{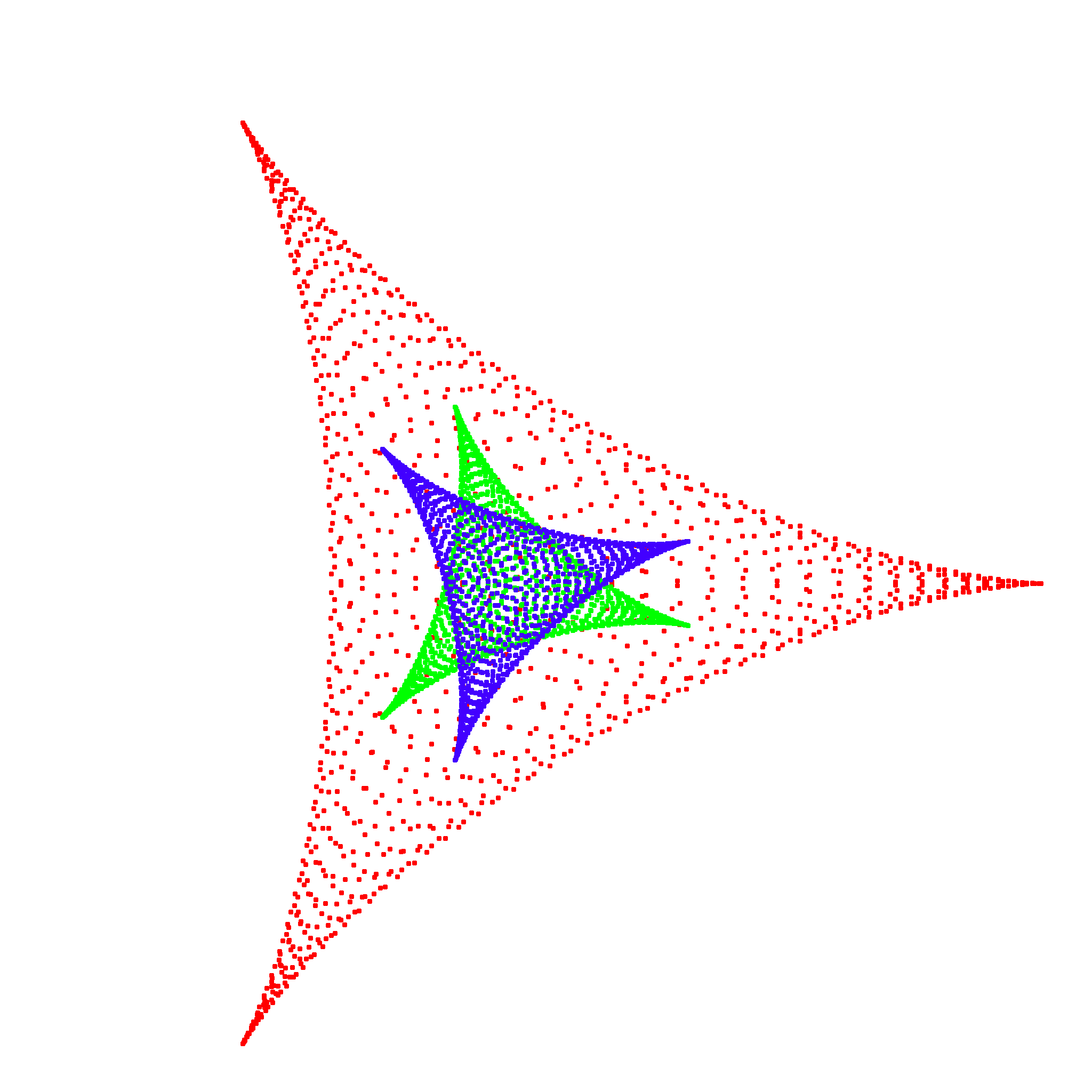}
        	\subcaption{$n=37367$, $\omega = 608$, $c=11$} 
        	\label{Figure:Finale:a}
\end{minipage}
~ 
\begin{minipage}[b]{0.32\textwidth} 
	\includegraphics[width=\textwidth]{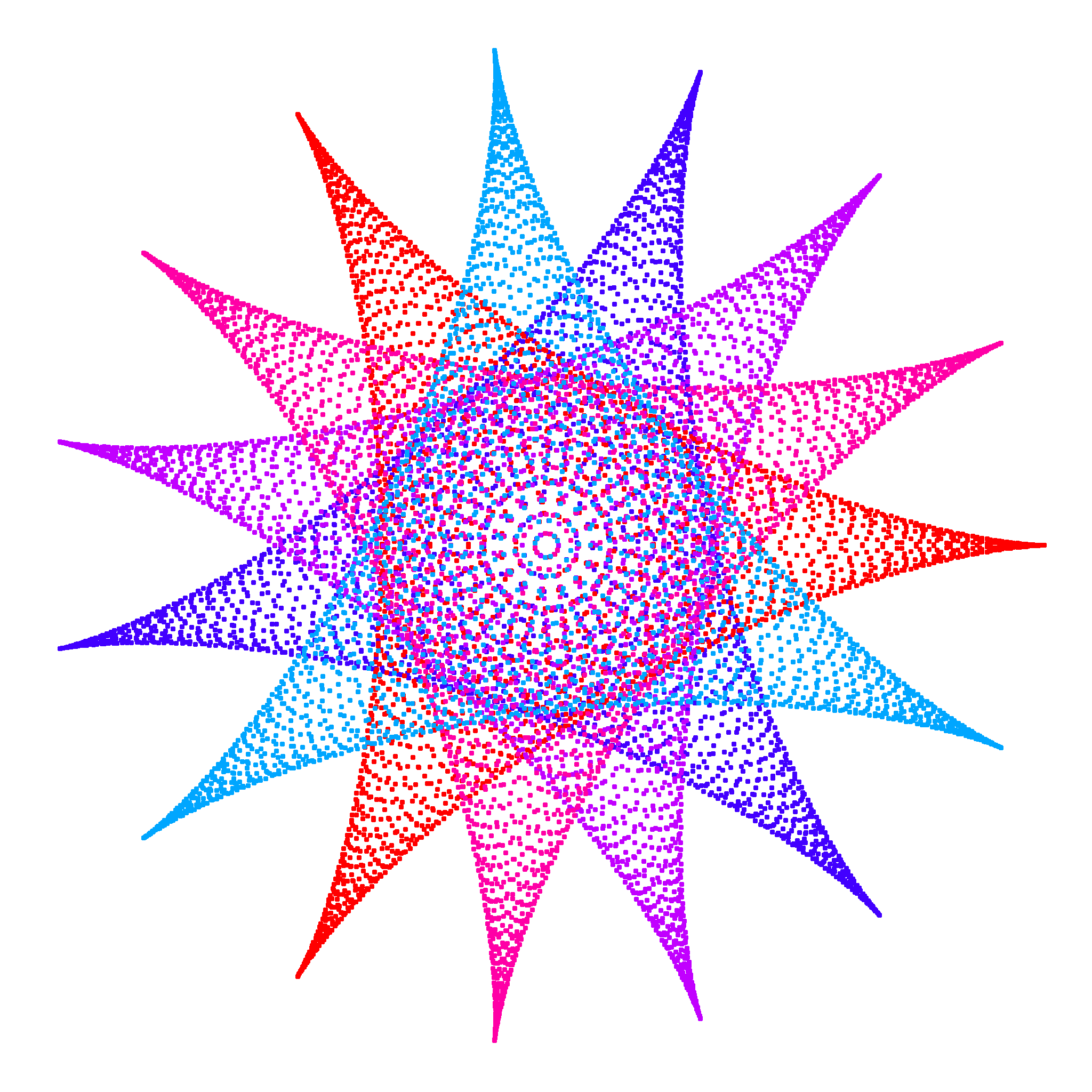}
        	\subcaption{$n=185925$, $\omega = 766$, $c=25$} 
        	\label{Figure:Finale:b}
\end{minipage}
~ 
\begin{minipage}[b]{0.32\textwidth} 
	\includegraphics[width=\textwidth]{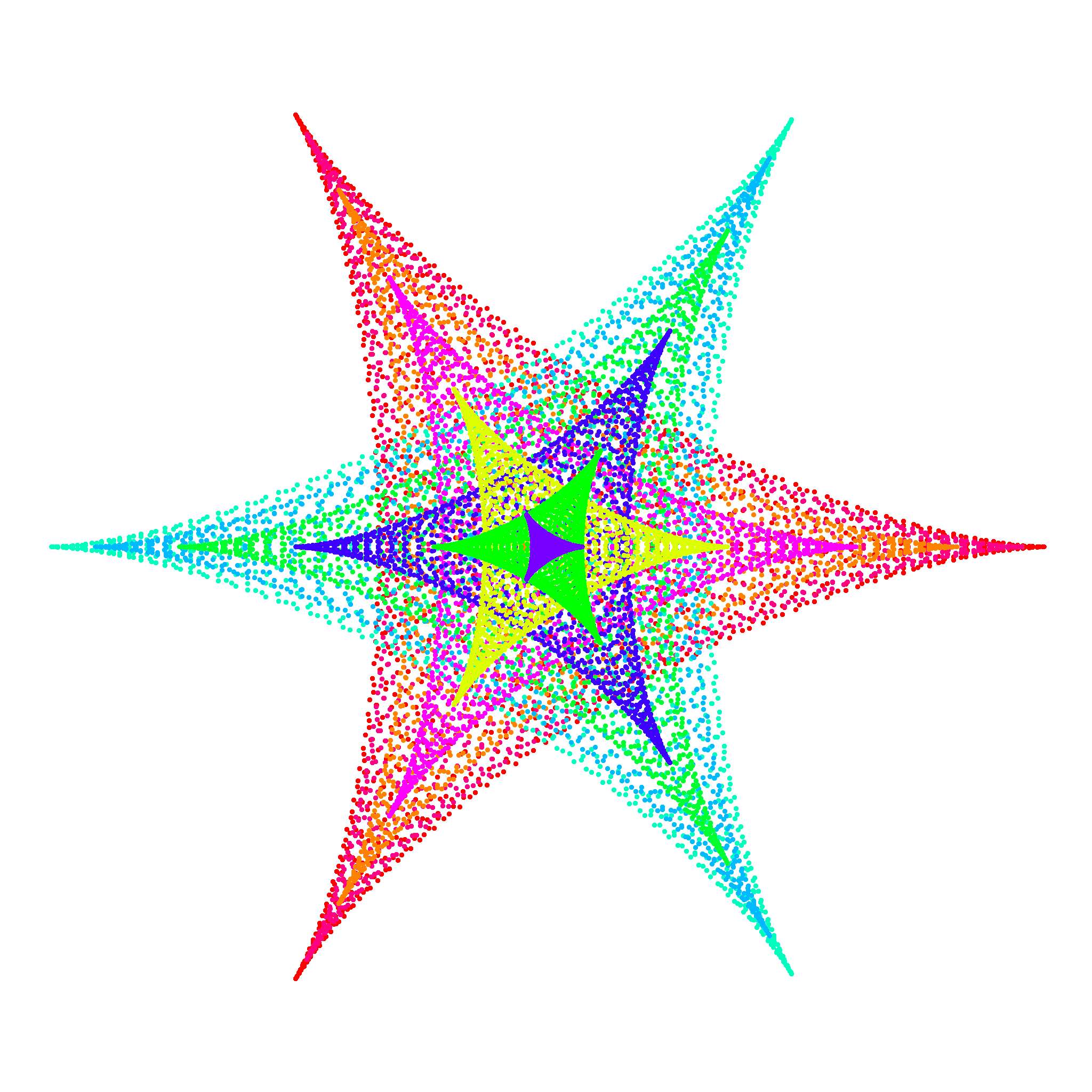}
        	\subcaption{$n=82677$, $\omega = 8147$, $c=21$} 
        	\label{Figure:Finale:c}		
\end{minipage}
\caption{Primitive graphical elements
(here, a ``filled'' deltoid) can form more elaborate plots.}
\label{Figure:Finale}
\end{figure}

With ellipses, hypocycloids, and so forth as primitive graphical elements, one can use the Chinese remainder
theorem to produce new images of startling complexity (as explained in \cite{Lutz}).  We content ourselves here
with a few aesthetically pleasing images produced in such a manner; see Figure \ref{Figure:Finale}.  
We encourage the reader to enjoy more examples by experimenting with our app \texttt{Gaussian Periods} (freely available at \url{http://www.elleneischen.com/gaussianperiods.html}).

\section*{Acknowledgments}
Ellen Eischen was partially funded by NSF grant DMS-1751281.  Stephan Ramon Garcia was partially funded by NSF grant DMS-1800123.  This material is based partly upon work supported by the NSF grant  DMS-1439786 and the Alfred P.~Sloan Foundation award G-2019-11406 while the authors were in residence at the Institute for Computational and Experimental Research in Mathematics (ICERM) in Providence, RI, during the Illustrating Mathematics program (Fall 2019).  We thank R.\ Lipshitz for help with the code.


\footnotesize
{\setlength{\baselineskip}{13pt} 
\raggedright				

} 

\end{document}